%% file: observationes.tex
\documentclass[12pt, letterpaper]{article}

\usepackage{amssymb}
\usepackage{booktabs}
\usepackage[T1]{fontenc}
\usepackage{graphicx}
\usepackage[latin,english]{babel}
\usepackage{ledmac,ledpar}
\lineation{section}
\usepackage{amsmath,graphicx,vmargin}
\setpapersize{USletter}

\newcommand{\exces}{\rotatebox{90}{\dag}}
\setmarginsrb{0.8in}{0.8in}{0.8in}{0.9in}{0pt}{0pt}{0pt}{24pt}

\begin{document}

\begin{center}
\textbf{\begin{Large}\textit{Observationes Cyclometric\ae{}} \\by Adam Adamandy Kocha\'nski -- Latin text with annotated English translation
\end{Large}}
\end{center}
\begin{center}
translated by Henryk Fuk\'s\\
Department of Mathematics,
Brock University, St. Catharines, ON, Canada\\
email \texttt{hfuks@brocku.ca}\end{center}

\noindent \textit{Translator's note:} The Latin text of \textit{Observationes} presented here closely  follows the original text published in \textit{Acta Eruditorum} 
\cite{Kochanski1685}. Punctuation, capitalization, and mathematical notation have been preserved. Several misprints which appeared in the original are also reproduced unchanged, but with a footnote indicating correction. Every effort has been made to preserve the layout of original tables.
The translation is as faithful as possible, often literal, and it is mainly intended to be of help to those who wish to study the
original  Latin text.

\setlength{\Lcolwidth}{0.465\textwidth}
\setlength{\Rcolwidth}{0.495\textwidth}

\begin{pairs}
\begin{Leftside}
\beginnumbering
\selectlanguage{latin}
\noindent\setline{1}
\input{part1-lat.tex}
%\endnumbering
\end{Leftside}

\begin{Rightside}
\beginnumbering
\selectlanguage{english}
\noindent\setline{1}
\input{part1-eng.tex}
%\endnumbering
\end{Rightside}

\Columns

\end{pairs}

\begin{center}
\textit{DIAMTERI AD PERIPHERIAM  CIRCULI}\\
Rationes Arithmetic\ae{}\footnote{\textit{Arithmetic ratios of diameter and circumference of a circle.} Ratios representing lower (``defective'') bounds are on the left, upper (``excessive'') bounds on the right.}
\end{center}
\begin{center}                        
\begin{tabular}{c|l||c|l}
 &  \hspace{6mm} \begin{large}\textit{Defectiv\ae{}} \end{large} &  &  \hspace{6mm} \begin{large}\textit{Excessiv\ae{}.} \end{large}\\ \hline
A & \hspace{1cm} 1. ad 3. \exces & Aa & \hspace{1cm} 1. ad 4. --- \\ \cline{2-2} \cline{4-4} 
B & \hspace{1cm} 8. ad 25. \exces & Bb &  \hspace{1cm} 7. ad 22. ---\\  \cline{2-2} \cline{4-4} 
Z & \hspace{1cm} 1.... 15.... 3. & Zz &  \hspace{1cm} 1.... 16.... 3. \\ 
C &  \hspace{5mm} 106. ad 333. \exces & Cc &  \hspace{1cm} 113. ad 355. ---\\ \cline{2-2} \cline{4-4}
Y & \hspace{1cm} 1.... 4697.... 3. & Yy &  \hspace{1cm} 1.... 4698... 3. \\  
D & 530762. ad 1667438 \exces & Dd & 530875. ad 1667793. ---\\ \cline{2-4}
X & \hspace{1cm} 1.... 5448\footnotemark.... 3. & Xx &  \hspace{1cm} 1.... 5449\footnotemark.... 3. \raisebox{1.2em}{\,}\\ 
E & Diam. 2945 294501. & Ee & Diam. 2945 825376. \\ 
 & Periph. 9252 915567 \exces &   & Periph. 9254 583360. ---\\ \cline{2-2} \cline{4-4} 
V & \hspace{1cm} 1.... 14774.... 3. & Vv &  \hspace{1cm} 1.... 14775.... 3. \\ 
F & Dia. 43 521624 105025. & Ff & Dia. 43 524569 930401 \\ 
 & Per. 136 727214 560643  \exces &   & Per. 136 736469 144003 ---\\ \hline
\end{tabular}
\end{center}
%%%%%Sorry, this must be done by hand, otherwise does not show
\footnotetext[7]{A misprint in the original text, should be 5548.}
\footnotetext[8]{Another misprint, should be 5549.}
\begin{center}
\newpage
 \textit{\begin{large}Harum q\ae{}dam minoribus terminis,\end{large}\footnote{\textit{Of these [ratios], some expressed in reduced form.} Here, cc, d, and e are
reduced forms of respectively Cc, D, and E, e.g. $\frac{71}{22\frac{3}{5}}=\frac{355}{113}$.}        }
\end{center}
\begin{center}
% use packages: array
\begin{tabular}{r|lcl}
cc & 22$\frac{3}{5}$ & ad & 71 --- \\ 
d. & 265381. & ad & 833719 \exces \\ 
e. & 30685681. & ad & 96401910 \exces\\ \hline
\end{tabular}
\end{center}

%  \begin{center}
%      \includegraphics[scale=0.5]{fig1.pdf}
%    \end{center}

\begin{pairs}
\begin{Leftside}
%\beginnumbering
%\memorydump
\selectlanguage{latin}
%\noindent\setline{1}
\input{part2-lat.tex}
%\endnumbering
\end{Leftside}

\begin{Rightside}
%\beginnumbering
%\memorydump
\selectlanguage{english}
%\noindent\setline{1}
\input{part2-eng.tex}
%\endnumbering
\end{Rightside}

\Columns
\end{pairs}
%\newpage
%\begin{samepage}
 \begin{center}
 \textit{\begin{large}Examen Rationum Cyclometricarum.\end{large}\footnote{\textit{Examination of cyclometric ratios.}} }
\end{center}
\begin{center}
% use packages: array
\begin{tabular}{c|r|r|rrr|r}\hline
Diam. & 100000 & 00000 & 00000 & 00000 & 00000 & Archimedis \\ 
Periph. & 314159 & 26535 & 89793 & 23846 & 26433 & -- Ratio. \\ \cline{1-6}
B. & 312500 & 00000 &  &   &   &   \\ 
  & 16519 & 26535 & \multicolumn{3}{|l|}{Defectus\footnotemark} & \\ 
Bb. & 314285 & 71428 & \multicolumn{3}{|l|}{}  &\\ 
  & 126 & 44892 & \multicolumn{3}{|l|}{Excessus.} & \\ \cline{1-6}
C. & 314150 & 94339 & \multicolumn{1}{|r|}{62264}  &   &   &   \\ 
  & 8 & 32196 & \multicolumn{1}{|r|}{27592} & \multicolumn{2}{|l|}{Defectus}  &   \\ 
Cc. & 314159 & 29203 & \multicolumn{1}{|r|}{53982} & \multicolumn{2}{|l|}{ } &   \\ 
  &   & 2667 & \multicolumn{1}{|r|}{64189} & \multicolumn{2}{|l|}{Excessus} &   \\ \cline{1-6}
D. & 314159 & 26535 & \multicolumn{1}{|r|}{81077} & \multicolumn{1}{|r|}{77120}  &   \\ 
  &   &   & \multicolumn{1}{|r|}{8715} & \multicolumn{1}{|r|}{46725} & Defect. &   \\ 
Dd. & 314159 & 26536 & \multicolumn{1}{|r|}{37862} & \multicolumn{1}{|r|}{02024} &   &   \\ 
  &   &   & \multicolumn{1}{|r|}{48068} & \multicolumn{1}{|r|}{78178} & Excess. &   \\ \cline{1-6}
E. & 314159 & 26535 & \multicolumn{1}{|r|}{89787} & \multicolumn{1}{|r|}{82814} &   &   \\ 
  &   &   & \multicolumn{1}{|r|}{5} & \multicolumn{1}{|r|}{41031} & Defect. &   \\ 
Ee. & 314159 & 26535 & \multicolumn{1}{|r|}{89796} & \multicolumn{1}{|r|}{49172} &   &   \\ 
  &   &   & \multicolumn{1}{|r|}{3} & \multicolumn{1}{|r|}{25326} & Excess. &   \\ \cline{1-6}
F. & 314159 & 26535 & \multicolumn{1}{|r|}{89793} & \multicolumn{1}{|r|}{23833} & \multicolumn{1}{|l|}{89913} &   \\ 
  &   &   & \multicolumn{1}{|r|}{ } & \multicolumn{1}{|r|}{12} & \multicolumn{1}{|l|}{36520} & Defectus. \\ 
Ff. & 314159 & 26535 & \multicolumn{1}{|r|}{89793} & \multicolumn{1}{|r|}{23855} & \multicolumn{1}{|l|}{91866} &   \\ 
  &   &   & \multicolumn{1}{|r|}{ } & \multicolumn{1}{|r|}{9} & \multicolumn{1}{|l|}{65432} & Excessus.\\ \hline \hline
\end{tabular}
\end{center}
%\end{samepage}
\footnotetext{\emph{Defect}, that is, the value of $\pi-B$. Similarly, \textit{excess} (lat. excessus) is the value of $Bb-\pi$.}
% \begin{center}
%      \includegraphics[scale=0.5]{fig2.pdf}
%    \end{center}
\begin{pairs}
\begin{Leftside}
%\beginnumbering
%\memorydump
\selectlanguage{latin}
%\noindent\setline{1}
\input{part3-lat.tex}
%\endnumbering
\end{Leftside}

\begin{Rightside}
%\beginnumbering
%\memorydump
\selectlanguage{english}
%\noindent\setline{1}
\input{part3-eng.tex}
%\endnumbering
\end{Rightside}

\Columns
\end{pairs}

 \begin{center}
     \includegraphics[scale=1.1]{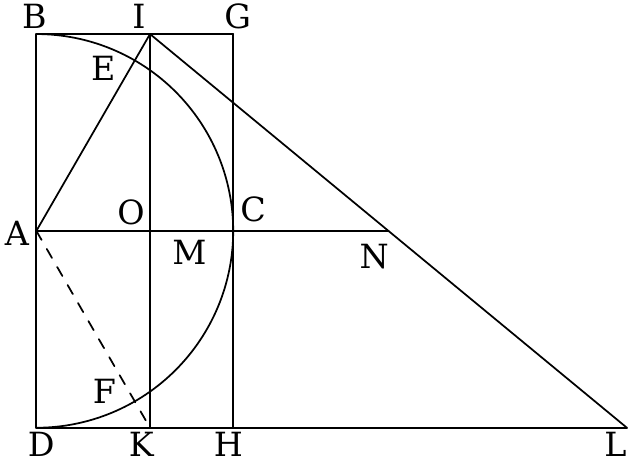}
   \end{center}

\begin{pairs}
\begin{Leftside}
%\beginnumbering
%\memorydump
\selectlanguage{latin}
%\noindent\setline{1}
\input{part4-lat.tex}
\endnumbering
\end{Leftside}

\begin{Rightside}
%\beginnumbering
%\memorydump
\selectlanguage{english}
%\noindent\setline{1}
\input{part4-eng.tex}
\endnumbering
\end{Rightside}

\Columns
\end{pairs}

\end{document}

%% file: part1-lat.tex
\pstart
\begin{center}
\emph{\begin{large}ADAMI ADAMANDI E SOCIET. JESU  \end{large}\\
Kochanski Dobrinniaci, Sereniss. Poloniarum Regis
Mathematici \&~Bibliothecari, OBSERVATIONES Cyclometric\ae{}, \\ad facilitandam Praxin accomodat\ae{};\\
ex Epistola ad Actorum Collectores.
}\end{center}
\mbox{\,}
\pend

\pstart
Qui Mathemata serio coluerit, nec tamen ad difficillima qu\ae{}que 
\& adhuc insoluta Problemata vires ingenii sui pertentandas 
censuerit, vix quenquam repertum esse existimo. Haud equidem 
diffiteor, me quoque olim eodem morbo laborasse, \& ut alia pr\ae{}teream,
in Circulo quidem quadrando, vel examinandis aliorum in eo 
conatibus, oper\ae{} non nihil collocasse. Non attinet hic enumerare
Methodos, quas ea in re secutus fueram: unam tantum, quam fortasse
quispiam felicius excolere poterit, commemorabo. Persuaseram mihi
conjectura quadam, possibiles esse aliquas Rectarum sectiones,
quarum segmenta invicem, \& cum aliis rectis Longitudine vel Potentia 
incommensurabilia essent, Circuli tamen Are\ae{}, vel Peripheri\ae{} partibus
Longitudine aut Potentia commensurarentur; ita ut inventa sectione
istiusmodi, liceret ex ea Tetragonismum expedire Geometrice, vel saltem
rationem Diametri ad Ambitum, in numeris ad lubitum maximis
supputare.
\pend

\pstart
Ad eam porro cogitationem videbar mihi non temere, sed illius
Quadratricis, a Dinostrato invent\ae{}, ductu devenisse. At cum ab
istis laboribus ad alia disparata studia animus avocaretur, illum tandem
adjeci, \& quidem magnorum Virorum exemplis incitatus, ad
investiganda compendia qu\ae{}dam Cyclometrica, Praxibus mechanicis
utilia, idque tam in Numeris, quam Lineis; quorum nonnulla hoc 
loco adferre lubet.
\pend

%% file: part1-eng.tex
\pstart
\begin{center}
\emph{
\begin{large}BY ADAM ADAMANY  FROM THE SOCIETY OF JESUS  \end{large} \\
Kocha\'nski of  Dobrzy\'n\footnote{Dobrzy\'n nad Wis\l\k{a} -- Kocha\'nski's birthplace, a town in Poland on the Vistula River, with settlement history dating back to 1065.},
 Mathematician and Librarian
of the Most Serene King\footnote{John III Sobieski (1629 -- 1696), from 1674 until his death King of Poland and Grand Duke of Lithuania.} of Poland,\\
Cyclometric OBSERVATIONS, \\accommodated
for easiness of practical use;
from a letter to fellow readers  of Acta.
}
\end{center}
\vskip 1cm
\pend

\pstart
I suppose one could hardy find anyone who would seriously cultivate knowledge\footnote{\textit{Mathemata} could mean both
knowledge or mathematics.} and who would nevertheless
not think that strengths of his talents are worth trying out on  difficult and yet unsolved problems.
For my part, I do not deny that I too  was once affected by the same weakness, and, to omit other things, 
I put not a small effort into squaring of a circle and in examination of works of others attempting it.
I does not belong here to list methods which I had followed in this matter: I will mention only one, which perhaps 
somebody luckier will be able to improve. I had convinced myself about a certain conjecture, namely that 
certain sections of a straight line are possible, whose fragments are incommensurable to each other and to other
straight lines in length and square, yet commensurable to parts of area or circumference in length or square;
so that by finding the section with this method, one might procure from it a quadrature of the circle  geometrically,
or at least compute the ratio of the diameter and circumference with as many digits as one likes.
\pend

\pstart
It seems that I have arrived to this idea not blindly, but guided by a
quadratrix\footnote{Quadratrix of Hippias is a curve with equation $y = x \cot ( \pi x/2 a)$. It can be used to solve the problem of
squaring the circle, although this is not a pure ``ruler and compass'' solution.}, 
invented by Dinostratus\footnote{Dinostratus (ca. 390 B.C. - ca. 320 B.C) was a Greek mathematician and geometer,
a disciple of Plato.}.
And while my mind was diverted from this work by other separate pursuits, eventually, 
inspired by examples of great men,  I turned to investigation of certain profits pertaining to
cyclometry, useful in mechanical practice, as much numerically as geometrically.
\pend

%% file: part2-lat.tex
\pstart
Methodicam pr\ae{}dictorum Numerorum Synthesin in \emph{Cogitatis,
\& Inventis Polymathematicis}, qu\ae{}, si DEUS vitam prorogaverit, utilitati
public\ae{} destinavi, plenius exponam; sufficiet interim ad eorum
notitiam insinuasse sequentia. Numeri Characteribus Z. Y. X. V.
tam simplicibus, quam geminatis insigniti, sunt \textit{Genitores}, e quorum
ductu, Numeri illis subjecti C. D. E. F. simplici, geminoque
charactere notati, procreantur hoc modo. Ratio 7. ad 22. Excessiva, 
ducta in Genitorem Z. 15; \& adjecto ad Productum Diametri, numero 1.
ad Peripheri\ae{} autem, hoc altero adjacente 3; constituit Rationem
C.106. ad 333, Defectivam: Genitor autem major Zz.16, ductus in eosdem
terminos Excedentes 7. ad 22, adjectisque ad horum Producta
numeris 1 \& 3, conficit Rationem CC. 113, ad 355. Excess:
\pend

\pstart

Similiter hi termini Excessivi 113, 355, multiplicati per Genitores
Y.Yy. videlicet 4697 \& 4698. servata adjectione numerorum 1. \& 3
ad Producta Diametri Peripheri\ae{}que, offerent terminos Rationum
D. \&Dd, qu\ae{} longe propius accedunt ad Archimedeam, a Ludolpho,
\& Gr\"umbergero nostro vastissimis expressam numeris. Eadem ratione
in reliquorum Terminorum genesi proceditur. Ut autem oculis
ipsis usurpare liceat, quantum exactitudinis adferant Rationes
ill\ae{}, visum est hoc loco adjicere Synopsin totius calculi, quo pr\ae{}dict\ae{}
Rationes ad Archimedeam, tanquam ad lapidem Lydium examinantur,
ut appareat, quantum sit uniuscujusque peccatum, defectu vel
excessu Peripheri\ae{} taxato in partibus Diametri totius, in particulas
Decimales subdivis\ae{}.
\pend

%% file: part2-eng.tex
\pstart
I will explain the aforementioned method more completely in \textit{Polymathic thoughts and inventions},
which work, if God prolongs my life, I have decided to put out for public benefit. In the meanwhile,
for acquaintance with this method, the following introduction will suffice.
Numbers denoted by both single and double characters Z, Y, X, V are \textit{Originators}, from which
numbers subjected to them, denoted by single and double characters A, B, C, D,  are derived this
way. Ratio of 7 to 22, excessive, is multipied by the Originator Z. 15, and with added  product of the diameter,
equal to 1, and 3,  close to circumference,  yields the defective ratio of 106 to 333.\footnote{%
Fraction $\frac{22}{7}$ is transformed into $\frac{22\cdot 15 +3}{7 \cdot 15 +1}=\frac{333}{106}$.}
Moreover, the major originator Zz.16, multiplied by the same exceeding bounds 7 and 22, and with  numbers 
1 and 3 added to the products, makes excessive ratio CC, 113 to 355.\footnote{This produces $\frac{22 \cdot 16 +3}{7 \cdot 16+1} = \frac{355}{113}$.}
\pend

\pstart
Similarly, those excessive bounds 113, 355, multiplied by originators Y and Yy, that is, 4697 and 4698, keeping addition of numbers 1 and 3, yield 
bounds\footnote{These bounds are $D=\frac{355 \cdot 4697+3}{113 \cdot 4697+1}=\frac{1667438}{530762}$ and $Dd=\frac{355 \cdot 4698+3}{113 \cdot 4698+1}=\frac{1667793}{530875}$.}
 on the ratio D and Dd, which come far closer to the Archimedean ratio, expressed by Ludolph\footnote{%
Ludolph van Ceulen (1540 -- 1610) was a German-Dutch  mathematician who calculated 35 digits of $\pi$.}
 and our Gru\"mberger\footnote{%
Christoph Grienberger  SJ (1561 -- 1636)  was an Austrian Jesuit astronomer, author of a catalog of fixed stars as well as optical and
mathematical works.}
 by great many digits. Remaining bounds are produced by proceeding in the same manner. In order to see  accuracy of these ratios with one's own eyes, it seemed fit to add in this place synopsis of all calculations, by which the aforementioned ratios are tested against Archimedean ratio like against the Lydian 
stone\footnote{Lydian stone (touchstone) -- stone used to test gold for purity.}, so that it would become evident how big was the error of each of these ratios, with defect or excess of the circumference expressed as parts of the entire diameter, and with digits divided into small groups.
\pend

%% file: part3-lat.tex
\pstart
Ex hac Tabella colligitur: imprimis quantitas Defectus, vel
Excessus cujusvis Rationis, taxata Fractione, cujus Denominator est
Diameter supremo loco posita, videlicet 1 cum tot zeris, quot
libet assumere: Numerator autem erit is, qui in eodem cum Denominatore gradu 
Decimali consistit. Sic Rationis C Defectum metitur
h\ae{}c Fractio $\frac{8}{100000}$ qu\ae{} exactior erit, si prolixior Denominator assumatur.
\pend

\pstart
%p 397
Colligitur ex eadem secundo. Rationes nobis exhibitas,
adeo compendiosas esse, ut earum nonnull\ae{}, duplo pluribus notis Archimedeis
\ae{}quivaleant ; quanquam ipsa Cc duplum earum excedat,
qu\ae{} proinde brevitate, nec non exactitudine sua, in Praxi c\ae{}teris 
pr\ae{}ferenda videatur, cui, dum quid accuratius qu\ae{}ritur, ipsa d succedat.
Pr\ae{}ter has quidem mihi suppetunt adhuc plures, consimili dote
pr\ae{}dit\ae{}, sed eas, ne nimius videar, alteri occasioni servandas existimo.
Concludam interim singulari quadam, \& ut ita dicam, curiosa Ratione,
qu\ae{} est 991 ad $3113\frac{991}{3113}$, qu\ae{} cum Archimedea consentit in octonis notis prioribus,
ac tum primum illam incipit excedere, minus quam 23 centesimis.
\pend

\pstart
\begin{center}
\textit{GRAMMIC\AE{} RATIONES CYCLOMETRIC\AE{},\\
Ad Usus Mechanicos.}
\end{center}
Harum quidem complures olim a me repert\ae{}; hoc tamen loco visum
mihi est eam tantum proponere, qu\ae{} huic Anno pr\ae{}senti,
quo ista scribimus, affinitate quadam conjuncta est.

Oporteat igitur Semiperipheri\ae{} B C D Rectam proxime \ae{}qualem reperire. Ducantur Tangentes B G, D H, quarum prior Radio
AC \ae{}qualis, \& jungantur GCH. Tum Radio CA secentur ex C arcus
utrinque \ae{}quales CE \& EF: quorum quivis complectetur Gradus 60,
reliqui autem BE, DF singuli gr. 30.  Agatur per E Secans AI, 
determinans Tangentem BI. Capiatur tandem HL, \ae{}qualis Diametro BD;
ac tum ducatur IL.
\pend

%% file: part3-eng.tex
\pstart
From this table one infers, first of all, quantities of the defect or excess of any ratio, estimated by a fraction whose denominator is the diameter placed in the initial position, with as many zeros as one wants to take. Numerator,  on the other hand, is this one, which takes the same position as the denominator. Thus the defect of the ratio C is estimated by the fraction $\frac{8}{100000}$, which would be more accurate if one took a longer denominator. 
\pend

\pstart
From the same table, a second thing is inferred. Ratios exhibited by us are so advantageous, that some of them are equivalent to
Archimedean ratio with twice as many digits as others; Yet Cc itself twice exceed others\footnote{That is, it exceeds B, Bb, and C in accuracy.}, hence  in practice by shortness and accuracy it seems to be preferred to others. If one sought a more accurate one, d would be a successor. Besides those,  I have  at hand indeed even more of them, similar in quality to the mentioned ones, but, in order not to appear excessive, I consider saving them for another occasion. I will, in the meanwhile, conclude with a certain singular, and so to speak curious ratio, which is 991 to 3113$\frac{991}{3131}$, which agrees with Archimedean in the first 8 digits, and then it starts to exceed it, by less than 23 hundredths.\footnote{%
Defining $r=\frac{3113\frac{991}{3131}}{991}=3.1415192677\ldots$,  we have $r-\pi \approx 0.23 \cdot 10^{-7}$.}
\pend

\pstart
\begin{center}
\textit{GEOMETRIC CYCLOMETRIC CONSTRUCTIONS,} \\
\textit{For Use by Mechanics.}
\end{center}
Of which several were once found by me. In this place, nevertheless, it seemed appropriate to present only one, associated with the current year, in which we write this.

It would be then required to find a straight line nearly equal to the semicircle BCD. Let tangent lines BG, DH be drawn, equal to the radius AC and connected by GCH. Then from C, let both parts of the arc be cut by CE and EF\footnote{This should likely be CF.}, equal to the radius CA. Each of them will embrace the angle of 60 degrees, while the remaining angles BE, DF will be 30 degrees each. Let a line AI be driven through E, determining the extent of the tangent BI. Finally, let HL be taken equal to the diameter BD; and then let IL be drawn.

\pend

%% file: part4-lat.tex
\pstart
Dico Inprimis IL \ae{}qualem esse
Semiperipheri\ae{} BCD proxime.
Demonstratur calculo Trigonometrico. Intelligatur
autem ducta esse IK, qu\ae{}
Tangentes BI, DK conjungat.\\
Quoniam ad Radium\\ 
\indent AB. \,\,\,\,100000 00000 00000.\\
Tangens gr. 30 est \\
\indent BI.   \,\,\,\,\,\,\,57735 02691 89626. Erit hujus\\
Compl. ad Radium, ipsa\\
\indent IG. \,\,\,\,\,\,\,42264 97308 10373. Igitur \\
Tota KH $\maltese$ HL, sive \\
\indent KL. \,\,2 42264 97308 10373. Ergo\\
IK q $\maltese$ XL q. \\ \indent 9 86923 17181 95572 75995 52843 99129.\\
Horum Radix est. \\
\indent IL. \,\, 3 14153 33387 05093. H\ae{}c autem\\
Deficit ab Archimedea. \\
\indent Z. \,\,\,\,\,\,\,\,\,\,\,\, .... 5 93148 84700.\\
Continetur in AB. vicib. \\
\indent X. \,\,  \,\,\,\,\,\,\,\,\,\,\,\,\,....\,\,\,\,\, ....\,\,\,\, 16859.
\pend

\pstart
Dico deinde, Peripheriam sic inventam, ab Archimedea ver\ae{} 
proxima deficere minori Ratione, ea, quam habet Unitas ad Decuplum
currentis Anni 1685 a Christo nato, \AE{}ra vulgari numerati;
majore autem, quam eadem habeat ad decuplum anni 1686,
proxime secuturi. Cum enim Peripheria nostra IL, ab Archimedea 
(pr\ae{}cedentis Tabul\ae{}) deficiat numero Z, qui totam Diametrum, numero
AB taxatam, metitur numero X: manifestum est, Unitatem ad numerum
pr\ae{}sentis Anni decuplum, videlicet 16850, majorem habere rationem
quam ad X, priore majorem: minorem autem, quam ad hunc
16860, per demonstrata Prop. 8 Quinti Elem. Euclid. E quibus
%akribeia = precision, exactness
%kai = but also
%eu)pori/a = easiness
Praxeos istius $\grave{\varepsilon} \upsilon \pi o \rho \acute{\iota}  \alpha$  $\kappa \alpha \grave{\iota}$ $\grave{\alpha} \kappa \varrho  \acute{\iota} \beta \varepsilon \iota \alpha$, cum intellectu comprehendi, 
tum memoria facile retineri poterit.
\pend

\pstart
Epimetri loco adjungam alteram Praxin Linearem Mechanicorum 
Circino opportunissimam, quod ea continua Diametri 
bisectione peragatur, sitque longe exactior pr\ae{}cedente: sic autem instituitur.
Dati Circuli Diametrum, Circino bisectioni destinato,
divide in partes 32. Talium enim Peripheria erit $100\frac{17}{32}$, hoc est, erit eorum
Ratio, 1024. ad 3217.  In Praxi igitur, Triplo Diametri,
sive partibus 96, adjiciend\ae{} erunt $\frac{4}{32}$, sive $\frac{8}{1}$ totius Diametri, \& insuper
Semis unius Trigesim\ae{} secund\ae{}, cum alterius Semissis particula
decima sexta. Hujus 
$\grave{\varepsilon} \gamma \chi \varepsilon \iota \rho \acute{\eta} \sigma \varepsilon \omega \varsigma$
%genitive of manipulation, operation
exactitudinem probat calculus,
quo provenit Peripheria P. 314160156.-- qu\ae{} Archimedeam
excedit numero Q .......891, qui minor est Defectu Z,
Peripheri\ae{} pr\ae{}cedentis. Quanquam nec istud subticendum sit,
istam Praxin in Majoribus Circulis potissimum locum habere, in
parvis oculorum effugere, quoad particulam, postremo addendam. 
\pend

%% file: part4-eng.tex
\pstart
I say in the first place that IL is nearly equal to the semicircle BCD. This is demonstrated by trigonometric calculations. Let us assume that
the line IK is drawn, which connects tangents BI, DK. 
%%%%%%%%%
 Since, with the radius\\ 
\indent AB. \,\,\,\,100000 00000 00000,\\
tangent of 30 degrees is \\
\indent BI.   \,\,\,\,\,\,\,57735 02691 89626. \\
Its complement\footnote{$1-\tan 30^{\circ}$.} to the radius, IG itself, will be\\
\indent IG. \,\,\,\,\,\,\,42264 97308 10373. Therefore, \\
together KH $\maltese$ HL, or \\
\indent KL. \,\,2 42264 97308 10373. Hence\\
IK q $\maltese$ XL\footnote{Obviously a misprint, should be KL instead of XL. IK q $\maltese$ KL q. means $IK^2 + KL^2 $.} q. \\ \indent 9 86923 17181 95572 75995 52843 99129,\\
of which the root is \\
\indent IL. \,\, 3 14153 33387 05093.\footnote{$IL=\sqrt{4+\left( 3- \frac{\sqrt{3}}{3}\right)^2}=\frac{1}{3}
\sqrt{120-18\sqrt{3}}$.} But this is\\
short of Archimedean by  \\
\indent Z. \,\,\,\,\,\,\,\,\,\,\,\, .... 5 93148 84700,\footnote{$Z=\pi-IL \approx 0.00005 93148 847$.}\\
contained in AB approx. \\
\indent X. \,\,  \,\,\,\,\,\,\,\,\,\,\,\,\,....\,\,\,\,\, ....\,\,\,\, 16859.\footnote{$X=\frac{1}{Z} \approx 16859$.}
%%%%%%%%%%%%%%
\pend

\pstart
I say then that the circumference found this way differs from the Archimedean ratio by less then the ratio of one to ten times the date of the current year 1685 after Christ, numbered with the common era, and more than one to ten times the date of the next year to come, 1686. When indeed our periphery IL is short of the Archimedean (of preceding table) by the number Z, which measures the whole diameter, expressed by the number AB, with the number X: it is clear that the ratio of the unity to ten times the current year, or 16850, is larger than the ratio of 1 to X, and less than 1 to 16860, as demonstrated by prop. 8 of the 5th book of Euclid\footnote{%
Prop. 8 of Euclid's book 5 says: \textit{Of unequal magnitudes, the greater has to the same a greater ratio than the less has; and the same has to the less
a greater ratio than it has to the greater.}\cite{Euclid}}.
From this exactness but also easiness of the construction, when embraced by the intellect, can easily be retained in memory.
\pend

\pstart
As a supplement, I will add another linear construction suitable for the compass of mechanics, which would be carried out by successive bisections of the diameter, and would be far more exact than the previous one: it is set up as follows. Given the diameter of the circle, intended for the bisection by compass, divide it into 32 parts.  Of such kind, the circumference will be $100\frac{17}{32}$, that is, will be the ratio of 1024 and 3217. In practice, therefore, to three diameters, or 96 parts,  $\frac{4}{32}$ will be added, or $\frac{8}{1}$ of the total diameter\footnote{Clearly, the author means $1/8$ here.}, and, moreover, one half of 32nd, with $\frac{1}{16}$ of another particle's
half.\footnote{This amounts to $\frac{96}{32}+\frac{4}{32}+\frac{1}{2} \cdot \frac{1}{32} + \frac{1}{32 \cdot 32}=\frac{3217}{1024}=3.1416015625$}
  Calculations prove exactness of this procedure, yielding Circumference P. 3.14160156, which exceeds Archimedean by the number Q .......891, which is smaller than the defect Z, of the preceding circumference\footnote{$Q=\frac{3217}{1024}-\pi \approx0.00000891$.}. In spite of this, one must not be silent about the fact that this construction has its place principally applied to larger circles, yet in smaller circles  it is beyond one's ability to see, especially with respect to the small particle added at the end.
\pend